\documentclass{amsart}

\usepackage{amssymb,amsmath,xcolor}
\usepackage{hyperref}

\newcommand\C{\mathbb C}

\newcommand\N{\mathbb N}
\newcommand\D{\mathbb D}
\newcommand\A{\mathcal A}
\newcommand\B{\mathcal B}
\newcommand\calF{\mathcal F}

\newcommand\Lip{\mathcal L}

\renewcommand\root{\frak{o}}

\newcommand\dist{\operatorname{d}}
\newcommand\im{\operatorname{ran}}
\newcommand\id{\operatorname{\mathbf 1}}
\newcommand{\cl}[1]{\overline{#1}}

\newcommand{\ran}[1]{\operatorname{ran}(#1)}
\def\f{\varphi}
\renewcommand\phi{\varphi}
\renewcommand\Re{\operatorname{Re}}
\newcommand\ds{\displaystyle}

\theoremstyle{plain}
\newtheorem{theorem}{Theorem}[section]
\newtheorem{proposition}[theorem]{Proposition}
\newtheorem{lemma}[theorem]{Lemma}
\newtheorem{corollary}[theorem]{Corollary}
\newtheorem{example}[theorem]{Example}

\theoremstyle{definition}
\newtheorem{defi}[theorem]{Definition}

\title[Composition operators on $\Lip_0$]{Composition operators on the little Lipschitz space of a rooted tree}
\author{Flavia Colonna}
\author{Rub\'en A. Mart\'inez-Avenda\~no}
\thanks{\mbox{}\\
  To Peter Rosenthal, in memoriam.
  \\ \\
  The second author is partially supported by the Asociaci\'on Mexicana de Cultura A.C}
\address{Department of Mathematical Sciences, George Mason University, Fairfax, VA, USA}
\address{Instituto Tecnol\'ogico Aut\'onomo de M\'exico, Mexico City, Mexico}
\email{fcolonna@gmu.edu}
\email{ruben.martinez.avendano@gmail.com}

\begin{document}

\maketitle

\begin{abstract}
  In this work, we study the composition operators on the little Lipschitz space $\Lip_0$ of a rooted tree $T$, defined as the subspace of the Lipschitz space $\Lip$ consisting of the complex-valued functions $f$ on $T$ such that $$\lim_{|v|\to\infty}|f(v)-f(v^-)|=0,$$
where $v^-$ is the vertex adjacent to the vertex $v$ in the path from the root to $v$ and $|v|$ denotes the number of edges from the root to $v$. Specifically, we give a complete characterization of the self-maps $\phi$ of $T$ for which the composition operator $C_\phi$ is bounded and we estimate its operator norm. In addition, we study the spectrum of $C_\f$ and the hypercyclicity of the operators $\lambda C_\phi$ for $\lambda \in \C$.
\end{abstract}

\section{Introduction}
In recent years the study of linear operators on discrete structures has emerged, partly inspired by the extensive study of such operators on many classical spaces in continuous environments. One very natural class of linear operators that has attracted the attention of researchers in the areas of complex analysis is the {\it composition operator}, that is, the linear operator $C_\f$ that maps a function $f$ defined on a set $D$  to the function $f\circ \f$, where $\f$ is a fixed self-map of $D$, called the {\it symbol} of the operator. 

In \cite{CoEa}, the authors studied the Lipschitz space $\Lip$ of an infinite rooted tree $T$ and in \cite{AlCoEa}, in joint work with Allen, the boundedness and the compactness of the composition operator acting on $\Lip$ were characterized. In addition, the isometries among the composition operators were identified and the spectrum of $C_\f$ was studied in detail when the operator is an isometry.

In this paper, we focus our attention on the composition operator acting on a closed separable subspace $\Lip_0$ of $\Lip$, known as the {\it little Lipschitz space} of the tree $T$.  

After presenting in Section~\ref{2} some preliminaries on trees and the definitions of the spaces $\Lip$ and $\Lip_0$, in Section~\ref{3}, we study the boundedness of the composition operator acting on $\Lip_0$ and provide estimates on the operator norm. In Section~\ref{4}, we provide several results on the spectra of $C_\f$. In Sections~\ref{5} and \ref{6}, we give several sufficient and some necessary conditions, respectively, on the hypercyclicity of a complex constant multiple of $C_\f$.

\section{Basic Definitions}\label{2}

A {\em tree} is a locally finite connected and simply-connected graph
$T$, which, as a set, we identify with the set of its vertices. Two
vertices $v$ and $w$ are said to be {\em neighbors} if there is an
edge connecting them, and we adopt the notation $v\sim w$. The
{\em degree} of a vertex is the number of its neighbors. A {\em path} between two vertices $v$ and $w$ is a finite sequence of distinct vertices $[v, v_1, v_2, \dots, v_{k-1}, w]$ such that $v\sim v_1\sim v_2  \sim \dots \sim v_{k-1}\sim w$. Since $T$ is a tree, the path between two vertices is unique. 

Given a tree $T$, we select a vertex $\root$ and designate it as the {\em root} of $T$. We call $T$ a  {\em rooted tree}. Given a vertex $w$ of such a tree $T$, a vertex $v$, with $v\neq w$, is called a {\em descendant} of $w$ if $w$ lies in the unique path from $\root$ to $v$. The vertex $w$ is called an {\em ancestor} of $v$. The {\em parent} of $v$ is the neighbor $v^-$ of $v$ which is an ancestor of $v$. The vertex $v$ is called a {\em child} of $v^-$.  The number $\dist(v,w)$ of edges of the path joining  two vertices $v$ and $w$ is called the {\em distance} between $v$ and $w$ and also the {\em
  length of the path} between $v$ and $w$. The length of a path joining the root $\root$ to a vertex $v$ is called the {\em length} of $v$ and is denoted by $|v|$.

Given a vertex $v$, the {\em sector determined by $v$} is the set
$S_v$ consisting of $v$ and all its descendants. Clearly,
$S_{\root}=T$. A {\em function on a tree} is a complex-valued function on the set of its vertices.

Throughout this article, we shall assume that $T$ is a rooted tree with root $\root$ whose vertices, save possibly $\root$, have degree greater than $1$ (and hence the tree is infinite). We shall adopt the notation $T^*$ for the set $T\setminus \{\root\}$.

The supremum norm of a bounded function $f$ on the tree is denoted by 
$\|f\|_\infty$.

\begin{defi}
We define the {\em Lipschitz space} of a tree $T$ to be the set  
$\Lip$ of functions $f: T \to \C$ such that
\[
\sup_{v \in T^*} |f(v)-f(v^-)| < \infty.
\]
\end{defi}

Colonna and Easley~\cite{CoEa} defined the Lipschitz space as the set of functions $f: T \to \C$ such that
\[
\sup_{v\neq w\in T} \frac{|f(v)-f(w)|}{\dist(v,w)} < \infty,
\]
In fact, they show that both definitions are equivalent and, furthermore,
\[
  \sup_{v \in T^*} |f(v)-f(v^-)| = \sup_{v \neq w\in T} \frac{|f(v)-f(w)|}{\dist(v,w)}
\]

Also, Colonna and Easley~\cite{CoEa} show that $\Lip$ is a Banach space when
endowed with the norm
\[
\| f \|:= |f(\root)| + \sup_{v \in T^*} |f(v)-f(v^-)|.
\]
In this paper, to make some calculations simpler, we endow the space
$\Lip$ with the equivalent norm
\[
\| f \|_\Lip:= \max\{ |f(\root)|,  \sup_{v \in T^*} |f(v)-f(v^-)|\}.
\]
With both norms, it can be shown that, for every $v, w \in  T$ we have
\begin{equation}\label{eq:norm_ineq}
|f(v)-f(w)| \leq \| f \|_\Lip \dist(v,w).
\end{equation}

Given a function $f: T \to \C$, we denote by $f'$ the function defined
as
  \[
  f'(v)=\begin{cases}
    f(v)-f(v^-), & \text{ if } v \neq \root, \text{ and }\\
    f(\root), & \text{ if } v = \root.
    \end{cases}
  \]
  Therefore,
  \[
    \| f \|_\Lip = \sup_{v \in T} |f'(v)|.
  \]

The space $\Lip_0$ is the subspace of $\Lip$ consisting of all
functions $f \in \Lip$ with
\[
\lim_{|v|\to \infty} |f(v)-f(v^-)|=0.
\]
It can be shown \cite{CoEa} that $\Lip_0$ is a closed subspace of
$\Lip$ and hence is itself a Banach space.

\begin{defi}
We say $\phi: T \to T$ is {\em Lipschitz} if
\[
\lambda_\phi:= \sup_{v\neq w} \frac{\dist(\phi(v),\phi(w))}{\dist(v,w)} < \infty.
\]
We say that $\lambda_\phi$ is the Lipschitz number of $\phi$.
\end{defi}

It can be shown~\cite{AlCoEa} that $\phi$ is Lipschitz if and only if 
\[
\sup_{v \in T^*} \dist(\phi(v),\phi(v^-)) < \infty.
\]
In this case, $\lambda_\phi=\sup_{v \in T^*}
\dist(\phi(v),\phi(v^-))$. Observe that $\lambda_\phi$ can be thought of as the maximum distance that the function $\phi$ takes apart two adjacent vertices.

We now prove an important property of $\Lip_0$ that will be useful in Sections \ref{4} and \ref{6}. (Also, compare with Theorem 2.3 in \cite{CoEa}.)

\begin{lemma}\label{le:finite_support}
The set $X$ of all functions in $\Lip_0$ with finite support is dense in $\Lip_0$.
\end{lemma}
\begin{proof}
Let $g \in \Lip_0$ and let $\epsilon >0$. We will show there is a
function $f \in X$ with $\| g -f\|_\Lip \leq \epsilon$.

Let $N \in \N$ be such that $|g'(v)| < \frac{\epsilon}{2}$ for every
vertex $v$ with $|v| \geq N$. Let $M \in \N$ such that
$\frac{|g(v)|}{M} < \frac{\epsilon}{2}$ for all vertices $v$ with
$|v|=N$ (this is possible, since there are only finitely many such
vertices). 

If $v$ is a vertex such that $|v|\ge N$, denote by $v_N$  the unique vertex of length $N$ on the path from $\root$  to $v$. We now define $f$ as follows:
\[
  f(v)=\begin{cases}
    g(v),& \quad\text{if }\ |v|\leq N,\\
\frac{M+N-|v|}{M} g(v_N), & \quad\text{if }\ N\leq |v|\leq N+M,\\
0,&\quad\text{if }\ |v|\geq N+M.\end{cases}
\]
Clearly the function $f$ is well defined and  belongs to $X$.

Since $\|g-f\|_\Lip=\| (g-f)'\|_\infty$, to complete the proof, it suffices to show that $|(g-f)'(v)|< \epsilon$ for all vertices $v \in T$. We distinguish the cases based on the length of $v\in T$. 

If $|v|=0$, then $v=\root$ and, since $(g-f)(\root)=0$, it follows that $(g-f)'(\root)=0$. 

If $0<|v|\leq N$, then $g(v)=f(v)$ and $g(v^-)=f(v^-)$ and hence $(g-f)'(v)=0$.

If $|v|=N+1$, then $f(v^-)=g(v^-)$ and $v_N=v^-$,  and hence 
\begin{eqnarray*}
|(g-f)'(v)| 
&=&|g(v)-f(v)| \\
&=&|g(v)-\frac{M-1}{M} g(v^-)| \\
&=& |g(v)-g(v^-)+\frac{1}{M} g(v^-) | \\
& \leq & |g(v)-g(v^-)|+\frac{1}{M} | g(v^-) | \\
& = & |g'(v)| + \frac{1}{M} |g(v^-)|\\
& < & \frac{\epsilon}{2} + \frac{\epsilon}{2} = \epsilon,
\end{eqnarray*}
since $|v|>N$ and by the choice of $M$.

If $N+1 < |v| \leq N+M$, then
\begin{eqnarray*}
|(g-f)'(v)|
&=& |g(v)-f(v) -g(v^-)+f(v^-)| \\
&<& |g(v)-g(v^-)| + |f(v)-f(v^-)| \\
& = & |g'(v)| + \left| \frac{M+N-|v|}{M} g(v_N)- \frac{M+N-|v^-|}{M} g(v_N) \right| \\
& = & |g'(v)| + \frac{1}{M} |g(v_N)|\\
& < & \frac{\epsilon}{2} + \frac{\epsilon}{2} = \epsilon,
\end{eqnarray*}
since $|v|>N$ and by the choice of $M$.

If $|v|>N+M$, then $f(v)=0=f(v^-)$ and hence
\[
|(g-f)'(v)| = |g'(v)| < \frac{\epsilon}{2}.
\]
 
We have shown that $|(g-f)'(v)| < \epsilon$ for all $v \in T$. Hence
$\| (g-f)' \|_\infty \leq \epsilon$ and thus $\| g-f \|_\Lip \leq \epsilon$ and the lemma follows.
\end{proof}

\section{Boundedness of composition operators on $\Lip_0$}\label{3}

Let $\calF$ be the set of all functions $f: T \to \C$. If  $\phi: T \to T$ is a function, we define the composition operator $C_\phi: \calF \to \calF$ as $C_\phi f = f \circ \phi$. This is clearly a linear transformation.

Allen, Colonna and Easley showed in \cite{AlCoEa},  that $C_\phi$ is bounded as an operator from $\Lip$ to $\Lip$ if and only if $\phi$ is a Lipschitz function. A natural question is to ask for necessary and sufficient conditions that guarantee that $C_\phi$ is bounded as an operator from $\Lip_0$ to $\Lip_0$. It is clear that if $C_\phi$ is bounded on $\Lip$ (equivalently, $\phi$ is Lipschitz), then $C_\phi$ will be bounded on $\Lip_0$ if and only if $C_\phi(\Lip_0) \subseteq \Lip_0$. In fact, using the Closed Graph Theorem it can be shown that $C_\phi$ is bounded on $\Lip_0$ if and only if $C_\phi(\Lip_0) \subseteq \Lip_0$, without using the boundedness of $C_\phi$ on $\Lip$.

  We start with a theorem that will be improved later, but which gives the taste of what is to come.

\begin{theorem}\label{theo:bounded}
Let $\phi: T \to T$ be a Lipschitz function. If
\begin{equation}\label{div} 
\lim_{|v| \to \infty} |\phi(v)|=\infty,
\end{equation}
then $C_\phi$ maps $\Lip_0$ into itself. 
\end{theorem}
\begin{proof}
Assume (\ref{div}) holds and let $f \in \Lip_0$. We need to show that $f \circ \phi \in \Lip_0$. Let $\epsilon >0$ and denote by $\lambda$ the Lipschitz number of
$\phi$. Observe that $\lambda>0$, since by assumption (\ref{div}), $\phi$ is not constant. Since $f
\in \Lip_0$, there exists $M \in \N$ such that
\[
|v|> M \quad \hbox{ implies } \quad | f(v) -f(v^-)| < \frac{\epsilon}{\lambda}.
\]
Now, since $\ds \lim_{|u| \to \infty} |\phi(u)|=\infty$, there exists $N \in \N$ such that 
\[
|u|> N \quad \hbox{ implies } \quad |\phi(u)|>M+\lambda.
\]

Let $u \in T$ such that $|u|> N$.  If $\phi(u)=\phi(u^-)$, then
trivially $ |f(\phi(u))-f(\phi(u^-))| < \epsilon$. On the other hand,
if $\phi(u)\neq \phi(u^-)$, there exists a unique path 
\[
[u_0, u_1, u_2, \dots, u_n]
\] 
with $u_0=\phi(u)$, $u_n=\phi(u^-)$ and $u_k\sim u_{k+1}$ for all $k=0,\dots,n-1$. 
 
Since $\phi$ is Lipschitz, we have that $n \leq \lambda$. Observe that, since $ |u_0| = |\phi(u)| > M + \lambda$, we must have that $|u_1|> M + \lambda -1$. But then, since $u_1 \sim u_2$, we must also have $|u_2|> M + \lambda - 2$. Continuing in this fashion we get that $|u_k|> M + \lambda -k$ for $k=0, 1, 2, \dots, n$. In particular, $|u_k|> M$ for $k=0, 1, 2, \dots, n$.

We then have that $|f(u_{k+1})-f(u_k)|< \frac{\epsilon}{\lambda}$ for
each $k=0, \dots, n-1$. But then,
\begin{eqnarray*}
|f(\phi(u))-f(\phi(u^-))| 
&\leq&  |f(u_0)-f(u_1)| + |f(u_1)-f(u_2)| + \dots \\
& & \quad \dots +  |f(u_{n-1})-f(u_n)| \\
&<& n \frac{\epsilon}{\lambda}  \\
&\leq& \epsilon.
\end{eqnarray*}
Thus, $f \circ \phi \in \Lip_0$, as desired.
\end{proof}

The following corollary gives easier-to-check conditions.

\begin{corollary}\label{cor_finite_inverse}
Let $\phi: T \to T$ be a Lipschitz function. If
$\phi^{-1}(\{v\})$ is finite for each $v \in T$, then $C_\phi$ maps
$\Lip_0$ into itself. In particular, if $\phi$ is injective, then
$C_\phi$ maps $\Lip_0$ into itself.
\end{corollary}
\begin{proof}
We only need to check that $\lim_{|v| \to \infty} |\phi(v)|=\infty$. Let $M>0$. Then, the set
\[
{\B} := \{  u \in T \, : \, | \phi(u) | \leq M \}=\bigcup_{|v|\leq M } \phi^{-1}(\{v\})
\]
is finite, by hypothesis. Choose $N>0$ such that for every $u \in
{\B}$, we have $|u|\leq N$. Consequently, if $|u|>N$, then $u \not\in {\B}$,
and hence $|\phi(u)|>M$, which is what we wanted.
\end{proof}

In Theorem~\ref{theo:bounded}, we should observe that we do not really need that $|\phi(v)|$ gets arbitrarily large for all vertices $v$ sufficiently far away from the root. For example, we may allow the function $\phi$ to be constant in a sector, while $|\phi(v)|$ gets arbitrarily large elsewhere, and we also obtain that $C_\phi$ maps $\Lip_0$ into itself. In fact, one can adapt the proof of Theorem~\ref{theo:bounded} to obtain the following result.

\begin{theorem}\label{LipsetA}
Let $\phi: T \to T$ be a Lipschitz function. Let
$\A=\{ v \in T \, : \, \phi \text{ is not constant on } S_v \}$. If
either $\A$ is finite or $\A$ is infinite and
\[
\lim_{\substack{|v|\to \infty \\ v \in \A}} |\phi(v)| = \infty,
\]
then $C_\phi$ maps $\Lip_0$ into itself. 
\end{theorem}
\begin{proof}
We wish to show that $f\circ\phi\in\Lip_0$ for each $f \in \Lip_0$. So let $f\in\Lip_0$ and fix $\epsilon >0$. If $\A$ is finite, then
there exists $M \in \N_0$ such that $\phi$ is constant at each $S_v$
with $|v|=M$. It then follows that for $|v|>M$, we will always have
$\phi(v)=\phi(v^-)$ and hence $|f (\phi(v))-f(\phi(v^-))|=0<\epsilon$,
which implies that $f \circ \phi \in \Lip_0$. Therefore we can now
assume that $\A$ is infinite.

The proof now is analogous to that of Theorem~\ref{theo:bounded}: Let $\lambda>0$ be the Lipschitz number of $\phi$. Since $f \in \Lip_0$, there exists $M
\in \N$ such that
\[
|v|> M \quad \hbox{ implies } \quad | f(v) -f(v^-)| < \frac{\epsilon}{\lambda}.
\]
Now, since $\ds \lim_{\substack{|v|\to \infty \\ v \in \A}} |\phi(v)|
= \infty$,  there exists $N \in \N$ such that
\[
u \in \A, \quad \text{ and } \quad  |u|> N \quad \text{ imply } \quad |\phi(u)|>M+\lambda.
\]

Let $u \in T$ such that $|u|> N+1$; hence, $|u^-|>N$. If $u^- \notin
\A$, then $\phi(u)=\phi(u^-)$ and then $ |f(\phi(u))-f(\phi(u^-))| =0
< \epsilon$ trivially.

Assume then that  $u^- \in \A$ and hence that
$|\phi(u^-)|>M+\lambda$. If $\phi(u)=\phi(u^-)$, then again,
trivially, $ |f(\phi(u))-f(\phi(u^-))|=0 < \epsilon$. On the other
hand, if $\phi(u)\neq \phi(u^-)$, 
the inequality $|f(\phi(u^-))-f(\phi(u))| < \epsilon$ can be obtained as it was done in the proof of Theorem~\ref{theo:bounded}.
Thus, $f \circ \phi \in \Lip_0$, as wanted.
\end{proof}

It turns out that the above condition is necessary as well.

\begin{theorem}\label{th:Cphi_not}
Let $\phi: T \to T$ be a Lipschitz function and let $\A$ be as in the statement of Theorem~\ref{LipsetA}.
If $\A$ is infinite and
\[
\lim_{\substack{|v|\to \infty \\ v \in \A}} |\phi(v)| \neq \infty
\]
then there exists $f \in \Lip_0$ such that  $f \circ \phi \notin \Lip_0$.
\end{theorem}
\begin{proof}

We first claim that there exists a vertex $w$ such that the set
$\phi^{-1}(\{w\}) \cap \A$ is infinite. Arguing by contradiction, assume that for every $w \in T$ the set $\phi^{-1}(\{w\}) \cap \A$ is finite. Let
$M>0$. Then the set
\[
{\B} := \{  u \in \A \, : \, | \phi(u) | \leq M \}=\bigcup_{|v|\leq M } \phi^{-1}(\{v\}) \cap \A
\]
is finite, by the assumption. Choose $N>0$ such that for every $u \in
{\B}$, we have $|u|\leq N$. Hence, if $u \in \A$ and  $|u|>N$, then $u
\not\in {\B}$, and hence $|\phi(u)|>M$, which shows that
\[ 
\lim_{\substack{|v|\to \infty \\ v \in \A}} |\phi(v)| = \infty,
\]
contradicting the hypothesis of the theorem. Hence,  there exists
a vertex $w$ such that the set $\phi^{-1}(\{w\}) \cap \A$ is infinite.

Now we set $f=\chi_{\{ w\}}$. Clearly, $f \in \Lip_0$. Let us show that $f \circ \phi \notin \Lip_0$. Indeed, since  $\phi^{-1}(\{w\})\cap \A$ is infinite, for each $M \in \N$ there exists a vertex $v \in \phi^{-1}(\{w\}) \cap \A$, with $|v|\geq M$. Since $\phi$ is not constant on $S_v$, there exists a vertex $u \in S_v$ such that $\phi(u)\neq w$. In fact, we may choose $u \in S_v$ such that $\phi(u)\neq w$ and $\phi(u^-)= w$. But then
  \[
    | (f\circ \phi)(u)-(f\circ \phi)(u^-) | = 1.
  \]
Hence, since for each $M \in \N$ we can find a vertex $u$ with $|u|\geq M$, such that $|(f\circ \phi)(u)-(f\circ \phi)(u^-)|=1$ it follows that $\lim_{|v|\to \infty} |(f\circ \phi)(v)-(f\circ \phi)(v^-)|\neq 0$ and hence $f \circ \phi \notin \Lip_0$.
\end{proof}

 We now show that the condition that $\phi$ is a Lipschitz function, which was used in all the previous results, is necessary as well. We thank a referee for providing us with the statement and the proof.

\begin{theorem}\label{th:CphiLip}
Let $\phi: T \to T$ be a function on the tree $T$. If $C_\phi$ is bounded on $\Lip_0$, then $C_\phi$ is bounded on $\Lip$. Hence $\phi$ is a Lipschitz function.  
\end{theorem}
\begin{proof}
  By the Closed Graph Theorem, to show the boundedness of $C_\phi$ on $\Lip$, it is enough to show that $C_\phi (\Lip) \subseteq \Lip$. Let $f \in \Lip$. Momentarily, fix a vertex $v \in T^*$ and let $N_v:=\max\{|\phi(v)|,|\phi(v^-)|\}$. Choose $m \in \N$ large enough so that
  \[
    \max\{ |f(w)| \, : \, |w|=N_v \} \leq m \, \| f \|.
    \]
We define a function $f_v \in \Lip_0$ in the following way:
  \[
    f_{v}(u)=\begin{cases}
               f(u), & \text{ if } |u|\leq N_v \\
               \frac{N_v+m-|u|}{m} f(u^*), & \text{ if } N_v \leq |u| \leq N_v + m,\\
               0, & \text{ if } N_v+m \leq |u|.
      \end{cases}
    \]
    where $u^*$ is the unique ancestor of $u$ with $|u^*|=N_v$. Observe that $|f_v(\root)|=|f(\root)|\leq \| f \|$; if $0< |u|\leq N_v$ then $|f_v(u) - f_v(u^-)|=|f(u)-f(u^-)| \leq \| f \|$; if $N_v < |u| \leq N_v+m$ then $|f_v(u) - f_v(u^-)|=\frac{|f(u^*)|}{m} \leq \| f \|$ by the choice of $m$; and if $|u|> N_v +m$ then $|f_v(u) - f_v(u^-)|=0$. Hence $\|f_v\| \leq \|f\|$. Furthermore, $f_v \in \Lip_0$.

Hence, since $|\phi(v)|\leq N_v$ and $|\phi(v^-)|\leq N_v$, we have
    \begin{align*}
      | f (\phi(v))-f(\phi(v^-))| & =|f_v(\phi(v))-f_v(\phi(v^-))| \\
                                 & =|(C_\phi f_v)(v)-(C_\phi f_v)(v^-)| \\
                                 & \leq \| C_\phi f_v \| \\
                                 & \leq \| C_\phi \|_{\Lip_0 \to \Lip_0} \,  \| f_v \| \\
                                 & \leq \| C_\phi \|_{\Lip_0 \to \Lip_0} \,  \| f \|,
    \end{align*}
    where the second inequality follows from the boundedness of $C_\phi$ on $\Lip_0$.
    
    Since we have obtained that $|f (\phi(v))-f(\phi(v^-))| \leq \| C_\phi \|_{\Lip_0 \to \Lip_0} \,  \| f \|$ for all $v \in T^*$, it follows that $f \circ \phi \in \Lip$, as desired. Hence $C_\phi$ is bounded on $\Lip$ and by \cite{AlCoEa}, we obtain that $\phi$ is a Lipschitz function.
  \end{proof}

  We can now combine the results above to obtain a full characterization.
  \begin{theorem}
    Let $\phi$ be a self map of the tree $T$ and let
    \[
      \A=\{ v \in T \, : \, \phi \text{ is not constant on } S_v \}.
    \]
    The composition operator $C_\phi$ is bounded on $\Lip_0$ if and only if 
    \begin{itemize}
\item $\phi$ is Lipschitz and
  \item either $\A$ is finite or $\A$ is infinite and $\ds
\lim_{\substack{|v|\to \infty \\ v \in \A}} |\phi(v)| = \infty$.
\end{itemize}
\end{theorem}
\begin{proof}
  If $C_\phi$ is bounded on $\Lip_0$, Theorem~\ref{th:CphiLip} implies that $\phi$ is Lipschitz. If $\A$ is finite, then we are done, so assume $\A$ is infinite. If
\[
\lim_{\substack{|v|\to \infty \\ v \in \A}} |\phi(v)| \neq \infty,
\]
then Theorem~\ref{th:Cphi_not} implies that $C_\phi$ is not bounded on $\Lip_0$, which contradicts the hypothesis.

On the other hand, if $\phi$ is Lipschitz and either $\A$ is finite or $\A$ is infinite and
\[
\lim_{\substack{|v|\to \infty \\ v \in \A}} |\phi(v)| = \infty\,
\]
then Theorem~\ref{LipsetA} implies that $C_\phi$ is bounded.
\end{proof}

We now obtain a sharp estimate for the norm of $C_\phi$ as an operator on $\Lip_0$. 

\begin{theorem}\label{th:bound}
Let $\phi$ be a self map of the tree $T$. If $C_\phi$
is bounded on $\Lip_0$, then
\[
  1+|\phi(\root)| \leq \| C_\phi \| \leq \max \{
  1+|\phi(\root)|, \lambda_\phi \},
\]
where $\lambda_\phi$ is the Lipschitz number of $\phi$.
\end{theorem}

\begin{proof}
Let $m=|\phi(\root)|$. Define the function $g:T \to \C$ as
\[
  g(v)=\begin{cases}
    1+|v|, & \text{ if } \ |v|\leq m, \\
    2m+1-|v|, & \text{ if } \ m \leq |v| \leq 2m+1, \\
    0, & \text{ if } \ 2m+1 \leq |v|.
  \end{cases}
\]
Then $g \in \Lip_0$, $\|g\|_\Lip=1$ and $(C_\phi
g)(\root)=g(\phi(\root))=1+|\phi(\root)|$. Hence, $\| C_\phi g \|_\Lip \geq
|g(\phi(\root))|=1+|\phi(\root)|$, which establishes the first inequality.

Now let $f \in \Lip_0$. The argument in \cite[Theorem 3.2]{AlCoEa} shows that
\[
  \sup_{v\in T^*} |(f\circ \phi)'(v)| \leq \lambda_\phi \sup_{v\in T^*} |f'(v)|, 
\]
and hence 
\begin{equation}\label{eq:norm1}
  \sup_{v\in T^*} |(f\circ \phi)'(v)| \leq \lambda_\phi \|f\|_\Lip.
\end{equation}
But by \cite[Lemma 3.1]{CoEa} (observe the norm in \cite{CoEa} is different, but the result still applies) we have
\[
  |f(\phi(\root))| \leq |f(\root)| + |\phi(\root)| \sup_{v\in T^*} |f'(v)|,
\]
and hence 
\begin{equation}\label{eq:norm2}
  |f(\phi(\root))| \leq \|f\|_\Lip + |\phi(\root)| \,  \| f \|_\Lip = (1+ |\phi(\root)|) \, \| f \|_\Lip.
\end{equation}
Combining equations \eqref{eq:norm1} and \eqref{eq:norm2} we obtain
\[
\| C_\phi f\|_\Lip \leq \max\{ 1+ |\phi(\root)|,\lambda_\phi \} \, \| f \|_\Lip,
\]
completing the proof.
\end{proof}

Observe that the proof of the above theorem also works on $\Lip$ and hence it gives an estimate for the norm of $C_\phi$ as an operator on $\Lip$.

Theorem~\ref{th:bound} allows us to obtain the exact value of the norm
when $\lambda_\phi \leq 1 + |\phi(\root)|$. In particular, if $\phi$
is constant, $\lambda_\phi=0$ and we obtain
\[
  \|C_\phi\|=1+|\phi(\root)|.
\]
This formula also holds if $\lambda_\phi=1$ (for example, if $\phi$ is an automorphism of $T$; i.e., $\phi$ is a bijective function of $T$ and $v \sim w$ if and only if $\phi(v) \sim \phi(w)$).
\smallskip

On the other hand, the following example shows that the case $\lambda_\phi > 1 + |\phi(\root)|$ may occur. 
\smallskip

Let $T=\N_0$ considered as a tree rooted at $0$ and let $\phi:T \to T$ be defined by
$\phi(0)=0$ and $\phi(n)=2n$ for $n \in \N$. Then, $\lambda_\phi=2$
but $|\phi(0)|=0$ and hence $\lambda_\phi > 1 + |\phi(0)|$.

In fact, in this case we can compute the norm of $C_\phi$ by
considering the function $f:\N_0 \to \C$ defined by 
\[
  f(n)=\begin{cases}
1&\quad\text{for } \ n=1,3,\\
2&\quad\text{for } \ n=2,\\
0&\quad\text{otherwise}. 
\end{cases}
\]
Then, $\|f\|_\Lip=1$ and
$(C_\phi f)(1)=2$ and $(C_\phi f)(n)=0$ if $n\neq 1$. Hence $\|C_\phi
f \|_\Lip=2$ and therefore $\|C_\phi\|=2=\lambda_\phi$. This also
shows that our estimate is sharp.

\section{The spectrum of $C_\f$ as an operator on $\Lip_0$}\label{4}

For the convenience of the reader, we recall the definitions of the spectra of a bounded operator $A$ on a Banach space $\B$ that we will use throughout this section.

The {\em spectrum}, denoted by $\sigma(A)$, is the set of complex
numbers $\lambda$ such that $A-\lambda$ is not invertible. The {\em
  point spectrum}, denoted by $\sigma_{\rm{p}}(A)$, is the set of complex
numbers $\lambda$ such that $(A-\lambda)f=0$ for some nonzero $f \in
\B$. The {\em approximate point spectrum}, denoted by
$\sigma_{\rm{ap}} (A)$, is the set of complex numbers $\lambda$ such
that there exists a sequence $\{ f_n \}$ of unit vectors in $\B$
with $(A-\lambda)f_n \to 0$, as $n \to \infty$. This is
equivalent to saying that $A-\lambda$ is not bounded
below (an operator $A$ is bounded below if
 there exists a constant $c>0$
such that $\| A f \| \geq c \| f \|$ for all $f \in
\B$). Lastly, the {\em compression spectrum}, denoted by
$\sigma_{\rm{c}}(A)$, is the set of complex numbers $\lambda$ such
that  $A - \lambda$ does not have dense range.

The following facts about the spectrum of an operator are well-known
(and they can be found in, for example, \cite[pp.~208, 349]{conway}
and in \cite{kubrusly}): 
\begin{itemize}
\item $\sigma(A)$ is a nonempty
compact subset of $\C$;
\item $\sigma_{\rm{p}}(A) \subseteq \sigma_{\rm{ap}}(A)$;
\item 
%the approximate point spectrum
 $\sigma_{\rm{ap}}(A)$ is a closed set containing the boundary of $\sigma(A)$;
 \item  $\sigma(A)=
\sigma_{\rm{ap}}(A) \cup \sigma_{\rm{c}}(A)$,
 but the union may fail to be
disjoint.
\end{itemize}

We now turn to investigate the point spectrum of a (bounded) composition
operator on $\Lip_0$. But first we make two general observations:

\begin{itemize}
\item For any $\phi : T \to T$, we have that $C_\phi \id = 1 \id$, where
$\id \in \Lip_0$ is the function with constant value $1$. Hence $1 \in
\sigma_{\rm{p}} (C_\phi)$, regardless of the choice of $\phi$.

\item If $\phi:T \to T$ is not surjective, then $0 \in \sigma_{\rm
  p}(C_\phi)$. Indeed, if  $w \notin \im(\phi)$, then $C_\phi
\chi_{\{w\}}=0$. (The converse is also true: if $C_\phi f=0$ for some
nonzero $f \in \Lip_0$, then $\phi$ is not surjective.)
\end{itemize}

In the next theorem we will find, under some conditions on $\phi$, other
numbers in the point spectrum of $C_\phi$. First, observe that  for any $\phi: T
\to T$, we have that $\chi_{\{ b\}} \circ \phi=\chi_{\phi^{-1}(\{b\})}$
for every $b \in T$: this fact will be useful in the next few
theorems. Also, observe that $\| \chi_{\{ b \}} \|_\Lip=1$ for every $b \in T$.

\begin{theorem}\label{th:point_spec}
Assume $\phi: T \to T$ is injective but not
surjective. If $C_\phi$ is bounded on $\Lip_0$, then $\{ 1 \} \cup \D
\subseteq \sigma_{\rm p}(C_\phi)$.
\end{theorem}
\begin{proof}
  That $1 \in \sigma_{\rm p}(C_\phi)$ has been showed already. Let
  $\lambda \in \D$ and assume $w \notin \im \phi$. Set
  \[
    f=\sum_{n=0}^\infty \lambda^n \chi_{\{\phi^n(w)\}}.
  \]
  Since $|\lambda|<1$ and $\| \chi_{\{\phi^n(w)\}} \|=1$, the sum above converges absolutely, and hence $f \in \Lip_0$.
  
  Now, if $n=0$, then $\chi_{\{\phi^n(w)\}} \circ \phi =\chi_{\{ w \}} \circ \phi=0$, since $w \notin \im \phi$.
  If $n>0$, then $\chi_{\{ \phi^n(w) \}} \circ \phi= \chi_{\{  \phi^{-1}(\{\phi^{n}(w)\}) \}} = \chi_{\{  \phi^{n-1}(w) \}}$,
  since $ \phi^{-1}(\{\phi^{n}(w)\}) = \{ \phi^{n-1}(w) \}$ for all $n \in \N$ by the injectivity of $\phi$. 

  Therefore,
  \[
    f \circ \phi = \sum_{n=0}^\infty \lambda^n \chi_{\{\phi^n(w)\}} \circ \phi  = \sum_{n=1}^\infty \lambda^n \chi_{\{\phi^{n-1}(w)\}} = \lambda \sum_{n=1}^\infty \lambda^{n-1} \chi_{\{\phi^{n-1}(w)\}}  = \lambda f,
  \]
  and hence, since clearly $f \neq 0$, it follows that $\lambda$ is an eigenvalue of $C_\phi$, with eigenvector $f$.
\end{proof}

In the proof of the theorem above, observe that, for each $\lambda \in
\D$, there are (at least) as many linearly independent eigenfunctions
as there are vertices $w$ not in the image of $\phi$. 

We now show a theorem establishing a bound on the point spectrum of
$C_\phi$. But first, we should point out that if $\phi:T \to T$ is
constant (say, $\phi(u)=w^*$ for every $u \in T$) then
$\sigma_{\rm p} (C_\phi)=\{0,1 \}$. Indeed, the observations
preceding Theorem~\ref{th:point_spec} show that $\{ 0, 1\} \subseteq
\sigma_{\rm p} (C_\phi)$. Now, assume that $\mu\neq 0, 1$ is an
eigenvalue of $C_\phi$ with eigenvector $g \in \Lip_0$. Then,
$g(\phi(u))=\mu g(u)$ for every $u \in T$. Evaluating at $w^*$ gives
$g(w^*)= \mu g(w^*)$. Since $\mu \neq 1$ it follows that
$g(w^*)=0$. But since $\mu\neq 0$, the equation $g(\phi(u))=\mu g(u)$
implies that $g(u)=0$ for every $u \in T$, and thus
$g=0$. Contradiction.

If $\phi$ is not constant, $\lambda_\phi \geq 1$, which gives a bound
for $\sigma_{\rm p}(C_\phi)$.

\begin{theorem}\label{specL0}
  Let $\phi: T \to T$ be a non constant Lipschitz function such
  that $C_\phi$ is bounded on $\Lip_0$. Then $\sigma_{\rm p}(C_\phi)$ is
contained in the closed disk centered at $0$ of radius $\lambda_\phi$.
\end{theorem}
\begin{proof}
Assume $\mu\in \sigma_{\rm p}(C_\phi)$ so that $C_\phi f = \mu f$ for some nonzero $f \in \Lip_0$. Then $f(\phi(u))=\mu f(u)$ for every $u \in T$. 

If $f$ is a (nonzero) constant function, then $\mu=1$ and since by assumption $\phi$ is not constant, it follows that $1 \leq \lambda_\phi$ and we are done in this case.

Next assume $f$ is nonconstant. Then 
\begin{equation}
  \sup_{u \in T^*} |f'(u)|>0.\label{sup>0}
\end{equation}

Let $u \in T^*$ and assume first that $\phi(u)\neq \phi(u^-)$. Then,
\begin{align*}
  |\mu (f(u)-f(u^-))|
  &= |f(\phi(u))- f(\phi(u^-))| \\
  &= \dist(\phi(u),\phi(u^-))\frac{|f(\phi(u)) -
    f(\phi(u^-))|}{\dist(\phi(u),\phi(u^-))}  \\
  &\leq \lambda_\phi \sup_{u \in T^*} |f'(u)|,
\end{align*}
where the last inequality follows from the definition of $\lambda_\phi$ and from \cite{CoEa}.

In the case when $\phi(u) = \phi(u^-)$, we have
\[
  |\mu (f(u)-f(u^-))| =| f(\phi(u))- f(\phi(u^-)) | = 0 \leq
  \lambda_\phi \, \sup_{u \in T^*} |f'(u)|.
\]
Hence, taking the supremum over all $u \in T^*$, we obtain
\[
  |\mu| \sup_{u \in T^*} |f'(u)|  =  |\mu| \sup_{u\in T^*} |f(u)-f(u^-)| \leq
 \lambda_\phi \sup_{u \in T^*} |f'(u)|.
\]
Hence, by (\ref{sup>0}), we obtain $|\mu| \leq \lambda_\phi$, as desired. 
\end{proof}

Note that Theorem~\ref{specL0} holds also for the point spectrum of $C_\phi$ as an operator on $\Lip$.

We now proceed to investigate the compression spectrum of $C_\phi$. 
\begin{theorem}\label{th:comp_spec}
Let $\phi: T \to T$. Then $C_\phi - \lambda$ has dense range if one of the following cases occurs:
\begin{enumerate}
\item $|\lambda|<1$ and $\phi$ is injective;
\item $|\lambda|>1$;
\item $\lambda\neq 0$ and for every $w \in T$, the set $\phi^{-n}(\{w\})$ is empty for some $n \in \N$.
\end{enumerate}
\end{theorem}
\begin{proof}
  Assume first that $|\lambda|<1$ and $\phi$ is injective. For each $w \in T$ we define
  \[
    f=\sum_{n=1}^\infty \lambda^{n-1} \chi_{\{\phi^n(w)\}}.
  \]
  Since $|\lambda|<1$, the above series converges absolutely and hence $f \in \Lip_0$.

  Now observe that, since $\phi$ is injective, we have $\chi_{\{\phi^n(w)\}} \circ \phi = \chi_{\{\phi^{n-1}(w)\}}$ for every $n \in \N$. Hence,
  \[
    C_\phi f = \sum_{n=1}^\infty \lambda^{n-1} \chi_{\{\phi^n(w)\}} \circ \phi = \sum_{n=1}^\infty \lambda^{n-1} \chi_{\{\phi^{n-1}(w)\}} = \chi_{\{w\}} + \sum_{n=1}^\infty \lambda^{n} \chi_{\{\phi^{n}(w)\}}.
  \]
  Since $\lambda f =  \sum_{n=1}^\infty \lambda^{n} \chi_{\{\phi^{n}(w)\}}$, it then follows that $(C_\phi - \lambda) f = \chi_{\{w\}}$. Hence $\chi_{\{w\}}$ is in the range of $C_\phi-\lambda$ and therefore, by linearity, all finitely supported functions in $\Lip_0$ are in the range of $C_\phi - \lambda$. By Lemma~\ref{le:finite_support}, the set of all finitely supported functions is dense in $\Lip_0$ and hence the range of $C_\phi-\lambda$ is dense in $\Lip_0$.\medskip

Next assume $|\lambda|>1$. For each $w \in T$ we define $f$ as
  \[
    f=\sum_{n=0}^\infty - \frac{1}{\lambda^{n+1}} \chi_{\phi^{-n}(\{w\})}.
  \]
Observe that $\phi^{-n}(\{w\})$ might be empty for some natural $n$, in which case the sum above is finite. If $\phi^{-n}(\{w\})$ is nonempty for every $n$, then the sum above converges absolutely, since $|\frac{1}{\lambda}|<1$ and $\| \chi_{\phi^{-n}(\{w\})} \|_\Lip \leq 1$ . In
either case, $f \in \Lip_0$. But again,
\[
C_\phi f = \sum_{n=0}^\infty - \frac{1}{\lambda^{n+1}}
  \chi_{\phi^{-n}(\{w\})} \circ \phi = \sum_{n=0}^\infty - \frac{1}{\lambda^{n+1}}
  \chi_{\phi^{-n-1}(\{w\})} = \sum_{n=1}^\infty - \frac{1}{\lambda^{n}}
  \chi_{\phi^{-n}(\{w\})}, 
\]
and
\[
  \lambda f = \sum_{n=0}^\infty - \frac{1}{\lambda^{n}} \chi_{\phi^{-n}(\{w\})}.
\]
Therefore,
\[
  (C_\phi - \lambda)f = \chi_{\{w\}}.
\]
As before, this implies that the range of $C_\phi-\lambda$ contains
all finitely supported functions in $\Lip_0$ and thus the range of
$C_\phi-\lambda$ is dense.  
\medskip

Finally, assume that $\lambda\neq 0$ and for every $w \in T$, the set $\phi^{-n}(\{w\})$ is empty for some $n \in \N_0$. For each $w \in T$ we define $f \in \Lip_0$ as
  \[
    f=\sum_{n=0}^{N-1} - \frac{1}{\lambda^{n+1}} \chi_{\phi^{-n}(\{w\})},
  \]
where $N \in \N$ is chosen so that $\phi^{-N}(\{w\})$ is empty. The computation above then yields that
\[
  (C_\phi - \lambda)f = \chi_{\{w\}},
\]
for every $w \in T$ and thus the range of $C_\phi-\lambda$ is dense.
\end{proof}

The theorem above shows the following result.

\begin{corollary}\label{cor:spec}
  Let $\phi:T \to T$ and assume $C_\phi$ is bounded on $\Lip_0$. Then
  \begin{enumerate}
    \item $\sigma_{\rm{c}}(C_\phi) \subseteq \cl{\D}$.
    \item If $\phi$ is injective, then $\sigma_{\rm{c}}(C_\phi)
      \subseteq S^1$.\label{circle}
    \item If $\phi$ is injective, then
      $\sigma(C_\phi)=\sigma_{\rm{ap}}(C_\phi)$.
    \end{enumerate}  
  \end{corollary}
  \begin{proof}
   The first two assertions follow trivially from parts (1) and (2) of Theorem~\ref{th:comp_spec} above. The third assertion follows from (\ref{circle}) and the fact that the approximate point spectrum of an operator is a closed set containing the boundary of the spectrum.
  \end{proof}

Observe that Theorem \ref{th:comp_spec} implies that  if $\phi$ is injective and for every $w \in T$ the set $\phi^{-n}(\{w\})$ is empty for some $n \in \N_0$, then $\sigma_{\rm{c}}(C_\phi)$ is empty. However, if $\phi$ is not injective, the compression spectrum is not empty and, furthermore, it is not contained in the unit circle, as the following proposition shows.
  
  \begin{proposition}\label{prop:not_dense_range}
If $\phi: T \to T $ is not injective,
    then $C_\phi$ does not have dense range. In other words, $0 \in
    \sigma_{\rm c}(C_\phi)$.
  \end{proposition}
  \begin{proof}
    Assume that $\phi(v)=\phi(w)$ for
    some $v\neq w$. Observe that if $g \in \im C_\phi$, then $g=f\circ
    \phi$ for some $f \in \Lip_0$. But then $g(v)=g(w)$. But the set
    $\{ g \in \Lip_0 \, : \, g(v)=g(w) \}$ is not dense in $\Lip_0$:
    for example, the function $\chi_{\{w\}}$ cannot be approximated by
    functions in this set. Hence $\im C_\phi$ cannot be dense.
  \end{proof}

 We now compute the spectrum for a concrete example.

\begin{example}\label{ex:example1}
Consider the tree $T=\N_0$, with $0$ being the root,
where two vertices $n$ and $m$ are adjacent if and only if $|n-m|=1$. If $\phi:T \to T$ is defined by $\phi(m)=2m+1$, then
  $\sigma(C_\phi)$ equals the closed disk, centered at the origin, of
  radius $2$.
\end{example}
\begin{proof}
  Since for any $n \in \N$ we have $\phi^n(m)=2^nm+ 2^n-1$, it follows
  that  $\lambda_{\phi^n}=2^n$ and $|\phi^n(0)|=2^n-1$. By Theorem
  \ref{th:bound}, it follows that $\| C_\phi^n
  \|=\|C_{\phi^n}\|=2^n$. Therefore, the spectral radius $r(C_\phi)$
  can be computed:
  \[
    r(C_\phi)=\lim_{n\to \infty} \| C_\phi^n \|^{1/n} = 2.
  \]
  We will now show that every point in the open disk centered at the
  origin of radius $2$ is an eigenvalue of $C_\phi$. 

  As observed in the observation preceding
  Theorem~\ref{th:point_spec}, since $\phi$ is not surjective we have
  that $0$ is an eigenvalue of $C_\phi$. Now, let
  $0<|\lambda|<2$. Choose $\mu \in \C$ so that $2^\mu=\lambda$ and
  define $f(m)=(m+1)^\mu:=\exp(\mu \log(m+1))$, where $\log$ is
  the natural logarithm.

  We first observe that $C_\phi f=\lambda f$. Indeed, for every $m \in \N$ we have
  \begin{align*}
    f(\phi(m))
    &=f(2m+1) \\
    &=(2m+2)^\mu \\
    &= 2^\mu (m+1)^\mu\\
    &= \lambda f(m).
  \end{align*}

  Now, we need to show that $f \in \Lip_0$. Fix $m \in \N$ and consider the
  function $g:[m,m+1] \to \C$ defined as
  \[
    g(t):=(t+1)^\mu:=\exp(\mu \log(t+1));
  \]
  (i.e., $g=f$ on $m$ and $m+1$). Then, by the Mean Value Inequality (e.g. \cite[Theorem 5.19]{rudin}), we obtain
  \[
    |f(m+1)-f(m)|= |g(m+1)-g(m)|\leq \max_{a\in [m,m+1]} \left|\frac{dg}{dt}(a)\right|.
  \]
  But for every $a \in [m,m+1]$, we have $\frac{dg}{dt}(a)=\mu (t+1)^{\mu-1}$, and therefore
  \begin{align*}
    \left|\frac{dg}{dt}(a)\right|
    &= |\mu| \, \left|(a+1)^{\mu -1}\right| \\
    &= |\mu| \, \left| \exp((\mu-1)\log(a+1)) \right| \\
    &= |\mu| \, \exp(\Re(\mu-1)\log(a+1)).
  \end{align*}
  But $\Re(\mu)=\frac{\log|\lambda|}{\log(2)}$ and hence
    \begin{align*}
    \left|\frac{dg}{dt}(a)\right|
      &= |\mu| \, \exp(\Re(\mu-1)\log(a+1))\\
      &= |\mu| \,
        \exp\left(\left(\frac{\log|\lambda|}{\log(2)}-1\right)
        \log(a+1)\right)\\
      &\leq  |\mu| \,
        \exp\left(\left(\frac{\log|\lambda|}{\log(2)}-1\right)
        \log(m+1)\right),
  \end{align*}
  since $\frac{\log|\lambda|}{\log(2)} - 1 < 0$ and $a\geq m$. Therefore,
  \[
    |g(m+1)-g(m)|\leq \max_{a\in [m,m+1]}
    \left|\frac{dg}{dt}(a)\right| \leq |\mu|  \exp\left(\left(\frac{\log|\lambda|}{\log(2)}-1\right)
        \log(m+1)\right).
  \]
  But the right-hand side of the above expression goes to $0$ as $m\to
  \infty$ (again, since $\frac{\log|\lambda|}{\log(2)} - 1 < 0$). Hence
   \[
    |f'(m+1)|=|f(m+1)-f(m)|  \to 0 \quad \text{ as } m \to \infty,
  \]
  and thus $f \in \Lip_0$. This completes the proof.
\end{proof}

Notice that in the proof of Example~\ref{ex:example1} we showed that $2 \D \subseteq \sigma_{\rm{p}}(C_\phi)$. Also, it turns out that for each $w \in T$, the set $\phi^{-n}(\{w\})$ is empty for some $n \in \N$ and that $\phi$ is injective. Hence, by Theorem~\ref{th:comp_spec} and Corollary~\ref{cor:spec},  $\sigma(C_\phi)=\sigma_{\rm{ap}}(C_\phi)=2 \cl{\D}$ and $\sigma_{\rm{c}}(C_\phi)$ is empty. 

\section{Hypercyclicity: necessary conditions}\label{5}

In this section, we investigate the hypercyclicity of multiples of the composition operator $C_\phi$ on $\Lip_0$. For the definition and basic facts about hypercyclicity, we refer the reader to the book \cite{GrEr-PeMa}.

First we notice that, as it is the case for composition operators on the Hardy-Hilbert space (e.g. \cite{shapiro}), if $\phi: T \to T$ has a fixed point, then $C_\phi$ cannot be hypercyclic. More is true.

\begin{proposition}\label{prop:periodic}
Let $\phi$ be a self map of $T$ such that
$C_\phi$ is bounded on $\Lip_0$. If $\lambda C_\phi$ is hypercyclic for some
$\lambda \in \C$, then $\phi$ does not have periodic points. In particular, $\phi$ does not have fixed points.
\end{proposition}
\begin{proof}
Suppose that there exists $v_0 \in T$ and $p \in \N$ such that
$\phi^p(v_0)=v_0$. Let $f \in \Lip_0$ be a hypercyclic vector for $\lambda C_\phi$. Then the
set
\[
\{  f, \lambda f \circ \phi, \lambda^2 f \circ \phi^2, \lambda^3 f \circ \phi^3, \dots  \}
\]
is dense in $\Lip_0$. Since $\Lip_0$ is a functional Banach space, it follows that
\[
\{  f(v_0), \lambda f (\phi(v_0)), \lambda^2 f (\phi^2(v_0)), \lambda^3 f(\phi^3(v_0)), \dots  \}
\]
is dense in $\C$. Since $\phi^p(v_0)=v_0$ this set equals
\[
\{  f(v_0), \lambda f(\phi(v_0)), \lambda^2 f(\phi^2(v_0)), \dots
\lambda^{p-1} f(\phi^{p-1}(v_0)), \lambda^p f(v_0),  \lambda^{p+1}
f(\phi(v_0)), \ldots \}.
\]
But this set equals
\[
  \bigcup_{j=0}^{p-1}\{ \lambda^{p k+j} f(\phi^{j}(v_0)) \, : \, k \in \N_0 \}
  \]
which clearly is not dense in $\C$. This is a contradiction, and hence $\phi$ does not have periodic points.
\end{proof}

As a consequence of the above proposition, we obtain that, in fact, the orbit of a vertex under $\phi$ moves away from the root.

\begin{corollary}\label{cor:phin_infty}
  Let $\phi: T \to T$ be such that $C_\phi$ is a
  bounded on $\Lip_0$. If $\lambda C_\phi$ is hypercyclic for some $\lambda \in \C$, then for every $w \in T$, we have
  \[
    \lim_{n\to \infty} |\phi^n(w)|=\infty.
  \]
\end{corollary}
\begin{proof}
  To show that $|\phi^n(w)| \to \infty$ we act by contradiction: assume there is a
  bounded subsequence $( \phi^{n_k}(w) )$. That is, there exists
  $M>0$ such that $|\phi^{n_k}(w)| < M$ for all $k \in \N$. Since there are only
  finitely many vertices at distance at most $M$ from the root, this implies
  that there exists $k\neq j$ such that $\phi^{n_k}(w) =
  \phi^{n_j}(w)$. But this implies that $\phi$ has a periodic point,
  which cannot occur since $\lambda C_\phi$ is hypercyclic. Therefore there
  cannot be bounded subsequences, and hence
  \[
    \lim_{n\to \infty} |\phi^n(w)|=\infty.
  \]
  as desired.        
\end{proof}

If for some different vertices $v_1$ and $v_2$ we have that $\phi(v_1)=\phi(v_2)$, then it is clear that $(C_\phi^n f)(v_1)= (C_\phi^n f)(v_2)$ for all $n \in \N$ and hence $C_\phi$ cannot be hypercyclic, since the orbit of $C_\phi$ cannot approximate functions which differ at $v_1$ and $v_2$. In fact, we had already shown this, as the proof of the next proposition makes clear.

\begin{proposition}\label{prop:not_injective}
Let $\phi$ be a noninjective self map of $T$ such
 that $C_\phi$ is bounded and let $\lambda \in \C$. Then $\lambda C_\phi$ is not hypercyclic.
\end{proposition}
\begin{proof}
Since  $\phi$ is not injective, then by Proposition~\ref{prop:not_dense_range}, $C_\phi$ does not have dense range, and hence neither does $\lambda C_\phi$. Therefore, $\lambda C_\phi$ cannot be hypercyclic.
\end{proof}

The previous proposition shows that we need to consider only composition operators induced by injective functions when studying hypercyclicity. The following result gives some information on the behavior of the orbits of every vertex and its parent under $\phi$, if $\lambda C_\phi$ is hypercyclic.

\begin{proposition}\label{prop:unbounded_distances}
Let $\phi$ be a self map of $T$ such that $C_\phi$ is bounded on $\Lip_0$. If $\lambda C_\phi$ is hypercyclic for some $\lambda \in \C$, then for all $w \in T^*$, the sequence $\left( \lambda^n \dist(\phi^n(w),
    \phi^n(w^-)) \right)$ is unbounded.
\end{proposition}
\begin{proof}
Assume that there exists a vertex $w$ such that
 $\left( \lambda^n \dist(\phi^n(w), \phi^n(w^-) )\right)$ is a bounded sequence; i.e. there exists $M>0$ such that
  \[
    \left| \lambda^n \dist(\phi^n(w), \phi^n(w^-)) \right| < M, 
  \]
  for all $n \in \N$. Since $\lambda C_\phi$ is hypercyclic, we may choose $f \in \Lip_0$ a
  hypercyclic vector with $\| f\|_\Lip < \frac{1}{2 M}$. Also, there exists $N\in \N$ such that
  \[
    \| \lambda^N C_\phi^N f - \chi_{\{w\}} \|_\Lip < \frac{1}{2}.
  \]
  But
  \begin{align*}
    \| \lambda^N C_\phi^N f - \chi_{\{w\}} \|_\Lip
    &\geq \left| (\lambda^N f \circ \phi^N- \chi_{\{w\}})'(w) \right|\\
    &=     \left| \lambda^N f(\phi^N(w))- \lambda^N f(\phi^N(w^-))-1\right|\\
    &\geq  1- |\lambda^N| \, \left| f(\phi^N(w))-f(\phi^N(w^-))\right|.
  \end{align*}
  The last two inequalities imply that
  \[
    |\lambda^N | \left| f(\phi^N(w))-f(\phi^N(w^-))\right| > \frac{1}{2}.
  \]
  But, by expression \eqref{eq:norm_ineq}, it follows that
  \[
   |\lambda^N | \left| f(\phi^N(w))-f(\phi^N(w^-))\right| \leq |\lambda^N | \, \| f \|_\Lip  \dist(\phi^N(w), \phi^N(w^-)) < \frac{1}{2M} M = \frac{1}{2},
  \]
  which is a contradiction.
  \end{proof}
  Observe that this proposition is only of interest in the case $|\lambda| \leq 1$. Indeed, observe that if $\lambda C_\phi$ is hypercyclic, then by the injectivity of $\phi$, we cannot have $\phi^n(w)=\phi^n(w^-)$ for any $n$ and hence  $\dist(\phi^n(w), \phi^n(w^-))\geq 1$. But then, if $|\lambda|>1$, we will have
  \[
    \left| \lambda^n \dist(\phi^n(w), \phi^n(w^-))\right| \geq |\lambda|^n,
    \]
    which is unbounded.

    As we mentioned before, the case $\lambda=1$ implies that no automorphism of $T$ induces a hypercyclic composition operator $C_\phi$. Compare this to the case of composition operators on the Hardy-Hilbert space \cite[Theorem 2.3]{BoSh}, in which there exist such operators.

A further necessary condition for the hypercyclicity of $\lambda C_\phi$ can be given. First, we show the following result.
      
    \begin{proposition}\label{last_prop}
If $\phi:T \to T$ does not have a periodic point, then for every
finite nonempty set $K \subseteq T$, there exists $N \in \N$ such that
$\phi^j(K) \cap K = \varnothing$ for $j \geq N$. In fact, if $|K|=n$,
then $\phi^j(K) \cap K \neq \varnothing$ for at most $n^2-n$ distinct
positive integers $j$.
\end{proposition}
\begin{proof}
Assume the conclusion is false and let $K$ be a nonempty finite set,
with $|K|=n$, such that $\phi^j(K) \cap K \neq \varnothing$ for at least
$n^2-n+1=n(n-1)+1$ distinct positive integers $j$. 

By the Pigeonhole Principle, since $K$ has $n$ elements, at least one
element of $K$, say $w$, has the property that there exist $n$
distinct positive integers $r_1, r_2, \dots, r_n$ and vertices $x_1,
x_2, \dots x_n \in K$ with $\phi^{r_j}(x_j)=w$. (If this didn't
happen, each of the vertices in $K$ would only be touched by a power
of $\phi$ at most $n-1$ times, for a total of $n(n-1)$ nonempty
intersections, which is impossible.)

Since $\phi$ does not have a periodic point, we have that $x_j \neq w$
for each $j$. Since $x_j \in K$ for each $j$ and $K$ has $n$ elements,
we must have that $x_s=x_t$ for some distinct indices $s$ and $t$. Denote 
$v:=x_s=x_t$. Then $\phi^{r_s} (v)=w$ and $\phi^{r_t}(v)=w$. Assume,
without loss of generality, that $r_s> r_t$. But then, if $m=r_s-r_t$
we have
\[
\phi^m (w)= \phi^m( \phi^{r_t} (v) )= \phi^{r_s} (v) =w
\] 
and hence $w$ is a periodic point. We have reached a contradiction.
\end{proof}

We should point out that the conclusion above is similar to the concept
of a ``run-away sequence'', introduced by Bernal-Gonz\'alez and
Montes-Rodr\'iguez in \cite{BeMo} in their study of the hypercyclicity of composition operators.

Combining Proposition~\ref{last_prop} with Proposition~\ref{prop:periodic}, we
obtain the following corollary.

\begin{corollary}
If $\lambda C_\phi$ is hypercyclic, then for every finite nonempty set
$K \subseteq T$, there exists $N \in \N$ such that $\phi^j(K) \cap K =
\varnothing$ for $j \geq N$.
\end{corollary}

\section{Hypercyclicity: some sufficient conditions}\label{6}

The previous results give necessary conditions for multiples of $C_\phi$ to be hypercyclic. In this section we give some sufficient conditions. We make use here of Lemma~\ref{le:finite_support} whenever we apply the Hypercyclicity Criterion. Recall that $X$ denotes the set of all complex-valued functions on $T$ with finite support.

\begin{theorem}\label{th:suff_hyp}
Let $\phi$ an injective self-map of a tree $T$ with no periodic points and suppose the operator $C_\phi$ is bounded on $\Lip_0$. If there exists an increasing sequence of nonnegative integers $(n_k)$ such that for every cofinite set $S \subseteq T$, there exists $N \in \N$ such that  $\phi^{n_k}(S) \subseteq S$ for every $k \geq N$, then $\lambda C_\phi$ is hypercyclic on $\Lip_0$ for each $|\lambda|>1$.
\end{theorem}
\begin{proof}
Let $A=\lambda C_\phi$. We use the Hypercyclicity Criterion, as described in \cite[p.~74]{GrEr-PeMa}.
\begin{enumerate}
\item Let $f \in X$ and let $M \in \N$ be such that $f(v)=0$ for all $|v|
\geq M$.  Consider the set $K:=\{ v \in T \, : \, |v| < M \}$, which
is finite, and let $S$ be its complement. By hypothesis, there exists $N_1 \in \N$ such that $\phi^{n_k}(S) \subseteq S$ for every $k \geq N_1$, and hence if $|v|\geq M$, then $|\phi^{n_k}(v)|\geq M$ for every $k \geq N_1$.

If $|v| < M$, since $\phi$ has no periodic points, there exists $R_v
\in \N$ such that $|\phi^k(v)|\geq M$ for all $k \geq R_v$. Since $K$
is finite, we can take $N_2:=\max\{R_v \, : \, v \in K \}$ and hence
$|\phi^k(v)|\geq M$ for all $k \geq N_2$ and all $v \in K$. In fact, since the sequence $(n_k)$ is increasing, we can guarantee that $|\phi^{n_k}(v)|\geq M$ for all $k \geq N_2$ and all $|v| < M$.

Hence, if $k \geq N:=\max\{N_1,N_2\}$, we have that $|\phi^{n_k}(v)|\geq M$ for all $v \in T$. It follows that $f(\phi^{n_k}(v))=0$ for all $k \geq N$. Hence $A^{n_k} f = \lambda^{n_k} C_{\phi^{n_k}} f =0$ if $k \geq N$,  and hence the first condition of the Hypercyclicity Criterion is satisfied.

\item Define a sequence of functions $B_n : X \to X$ as follows. For each $f \in X$ we set
  \[
 (B_n f)(v):=\begin{cases}
   \frac{1}{\lambda^n}f(\phi^{-n}(v)), & \text{ if } v \in \ran{\phi^n}, \\
   0, & \text{ if not.}
      \end{cases}
    \]
    Observe that, by injectivity, there is a unique value of $\phi^{-n}(v)$ if $ v \in \ran{\phi^n}$. 

    Since $f \in X$, then $f$ is bounded. This implies that
    \[
      |(B_n f)(v)|\leq \frac{1}{|\lambda|^n} \| f \|_\infty
    \]
    for all $v$, and hence $\| B_n f \|_\Lip \leq \frac{2}{|\lambda|^n} \| f \|_\infty$. Therefore $\| B_{n_k} f \| \to 0$, as $k \to \infty$, and the second condition of the Hypercyclicity Criterion is satisfied.
    \item Lastly, it is clear that $(A^n B_n f)(v)= \lambda^n (B_n f)(\phi^n(v))=f(v)$ and hence $A^{n_k} B_{n_k} f \to f$ as $k \to \infty$.
 \end{enumerate}
By the Hypercyclicity Criterion, $A$ is weakly mixing, and hence, hypercyclic.
\end{proof}

If the sequence $(n_k)$ is the sequence of all natural numbers, Kitai's Criterion implies that $C_\phi$ is actually mixing \cite[p.~71]{GrEr-PeMa}. The following proposition gives alternative ways in which the hypotheses in the theorem above hold.

\begin{proposition}\label{prop:equivalent}
Let $\phi$ be an injective self map of $T$ with no
periodic points. Then the following conditions are equivalent:
\begin{enumerate}
  \item For every cofinite set $S \subseteq T$, there exists $N \in
    \N$ such that $\phi^k(S) \subseteq S$ for every $k \geq N$.
  \item There exists an increasing sequence of nonnegative integers $(n_k)$ such that  for every cofinite set $S \subseteq T$, there exists $N \in \N$ such that $\phi^{n_k}(S) \subseteq S$ for every $k \geq N$.
  \item For every $v \in T$, the set $\{n \in \N \, : \, \phi^{-n}(\{ v \}) \neq \varnothing \}$ is finite.
\end{enumerate}
\end{proposition}
\begin{proof}
That (1) implies (2) is obvious: just take $n_k=k$.

Now assume (2). For clarity of notation, observe that since $\phi$ is injective, if $\phi^{-1}(\{v\})$ is nonempty, we can refer to the unique element in $\phi^{-1}(\{v\})$ as simply $\phi^{-1}(v)$. With this notation, suppose there is a vertex $v \in T$ such that the set $\{n \in \N \, : \, \phi^{-n}(\{v\})\neq \varnothing \}$ is infinite. Then the set $S:=T\setminus\{ \phi^{-1}(v) \}$  is cofinite and therefore there exists $N \in \N$ such that $\phi^{n_k}(S)\subseteq S$ for all $k \geq N$. By the assumption, $\phi^{-n_N-1}(\{v\})$ is nonempty and, furthermore $\phi^{-n_N-1}(v) \neq  \phi^{-1}(v)$, since $\phi$ has no periodic points. Therefore $\phi^{-n_N-1}(v) \in S$ and hence, $\phi^{-1}(v)=\phi^{n_N}(\phi^{-n_N-1}(v)) \in S$, which is a contradiction. Therefore (3) holds.

Now assume (3). Arguing by contradiction, assume that there exists a cofinite set $S \subseteq
T$ and an increasing sequence of natural numbers $\{n_j\}$ such that
$\phi^{n_j} (S) \not\subseteq S$. Since $S$ is cofinite, that means
there is a vertex $w \not\in S$ and an infinite number of indices $j$
such that there exists $s_j\in S$ with $\phi^{n_j}(s_j)=w$. Since
$\phi$ has no periodic points, this implies that the set $\{n \in \N \, : \, \phi^{-n}(\{w\})\neq \varnothing \}$ is infinite, which violates (3). Hence (1) holds.
\end{proof}

Combining Proposition~\ref{prop:equivalent}, Theorem~\ref{th:suff_hyp} and the comment after Theorem~\ref{th:suff_hyp}, we obtain the following corollary.

\begin{corollary}\label{cor:mixing}
  Let  $\phi$ be an injective self map of $T$ with no periodic points. If any (and hence all) of the conditions of Proposition~\ref{prop:equivalent} hold, then $\lambda C_\phi$ is mixing on $\Lip_0$ for each $|\lambda|>1$.
  \end{corollary}

  An alternative, and maybe easier way, in which we can check hypercyclicity for concrete composition operators, will be given by the following statement.

  \begin{proposition}\label{prop:nonreversible}
    Let $\phi$ be an injective self map of $T$ with no periodic points. If there exists $N \in \N$ such that $|v|< |\phi(v)|$ for every $v \in T$ with $|v| \geq N$, then $\{n \in \N \, : \, \phi^{-n}(\{v\}) \neq \varnothing\}$ is finite for every $v \in T$.
  \end{proposition}
  \begin{proof}
Assume that there is $v \in T$ such that $\{n \in \N \, : \, \phi^{-n}(\{v\}) \neq \varnothing\}$ is infinite. Since $\phi$ is injective, we refer to the single element in the set $\phi^{-n}(\{v\})$ simply as $\phi^{-n}(v)$. Since $\phi$ has no periodic points, it must be the case that $|\phi^{-n}(v)|< N$ for at most finitely many values of $n$. Hence, by renaming if necessary, we may assume that there is  $v \in T$ such that $\{n \in \N \, : \, \phi^{-n}(\{v\}) \neq \varnothing\}$ is infinite and $|\phi^{-n}(v)|\geq N$ for all $n \in \N$. But then the condition implies that
 \[
      |v| > |\phi^{-1}(v) | > |\phi^{-2}(v) | >|\phi^{-3}(v) | >|\phi^{-4}(v) | > \cdots,
    \]
    which is impossible, since all numbers above are nonnegative integers. This finishes the proof.    
  \end{proof}

  The implication above is strict. Indeed, let $T$ be the tree constructed as follows. Let the vertices of $T$ be the disjoint union
 \[
    \N_0 \cup \N^0\cup \N^2 \cup \N^4 \cup \N^6 \cup \cdots,
  \]
  where $\N_0$ is the set of nonnegative integers, and for each even $k \in \N_0$ we set $\N^k:=\{1_k, 2_k, 3_k, \cdots \}$  as a disjoint copy of the set of natural numbers, where we index them for identification purposes only. Two vertices $u$ and $v$ in $\N_0$ are adjacent in $\N_0$ if $|u-v|=1$, as usual; two vertices  $u_k$ and $v_k$ in $\N^k$ are adjacent in $\N^k$ if $|u-v|=1$; and we also set $1_k \in \N^k$ to be adjacent to $k \in \N_0$. We set the root to be $0 \in \N_0$. Define the function $\phi: T \to T$ as follows. If $v \in \N_0$, we set $\phi(v)=v-1$ if $v$ is odd, and $\phi(v)=1_v$, if $v$ is even. If $v_k \in \N^k$, we set $\phi(v_k)=(v+1)_k$. By construction, the function $\phi$ is injective and has no periodic points. Furthermore, the set $\{ n \in \N \, : \, \phi^{-n}(\{v\}) \neq \varnothing\}$ is finite for each $v \in T$, but $|\phi(v)| < |v|$ for each odd $v \in \N_0$. Hence the implication in Proposition~\ref{prop:nonreversible} is strict.
  
Nevertheless, we obtain the following corollary which gives a condition that is easy to check to show that some multiples of $C_\phi$ are mixing.

\begin{corollary}
  Let $\phi$ be an injective self map of $T$. If there exists $N \in \N_0$ such that $|u| < |\phi(u)|$ for every $u \in T$ with $|u|\geq N$, then $\lambda C_\phi$ is mixing on $\Lip_0$ for each $|\lambda|>1$.
  \end{corollary}

  For an easy application, consider Example \ref{ex:example1}: in that case, the function $\phi$ is injective and has no periodic points, and  $|u|<|\phi(u)|$ for every $u$. Hence $\lambda C_\phi$ is hypercyclic for each $|\lambda|>1$. 
  As we will see in Example~\ref{ex:example2}, the operator $C_\phi$ is in fact hypercyclic.

Incidentally, by Corollary~\ref{cor:mixing}, the operator $\lambda C_\phi$ is hypercyclic for each $|\lambda|>1$ also for the tree $T$ and function $\phi$ constructed after Proposition~\ref{prop:nonreversible}.  Nevertheless, in this case $C_\phi$ is not hypercyclic by Proposition~\ref{prop:unbounded_distances}.

  The result above can be extended, with additional conditions, for $|\lambda|\leq 1$. We thank an anonymous referee for providing us with the following result, which we simplified slightly.

  \begin{theorem}\label{th:hyper_referee}
    Let  $\phi$ be an injective self map of $T$ which has no periodic points and let $\lambda \in \C$ with $|\lambda|\leq 1$. Define, for each $v \in T$ and each $n \in \N$, the number
    \[
      m(n,v):=\min{\{\dist(\phi^n(u),\phi^n(v)) \, : \, u \in T, u \neq v \}}.
      \]
Assume that 
    \begin{enumerate}
      \item for all $v \in T$, we have $|\lambda|^n m(n,v) \to \infty$  as $n \to \infty$, and
\item there exists $c \in \N$ and $N \in \N$ such that for all $n \geq N$ and $v \in T$, we have 
  \[
    |\phi^n(v)| + c \geq m(n,v).
    \]
\end{enumerate}
Furthermore, if $|\lambda|=1$, assume that for each $v\in T$, the set $\{ n \in \N \, : \, \phi^{-n}(\{v\}) \neq \varnothing
\}$ is finite. Then $\lambda C_\phi$ is mixing on $\Lip_0$.
\end{theorem}
\begin{proof}
Again, we make use of the Hypercyclicity Criterion \cite[p.~74]{GrEr-PeMa}. Recall that $X$ is the set of all complex-valued functions on $T$ with finite support and $\|\chi_S\|_\Lip=0$ or $1$, depending on whether the set $S$ is empty or not.
\begin{enumerate}
\item For every $v \in T$ and every $n \in \N$, we have that  $(\lambda C_\phi)^n \chi_{\{v\}} = \lambda^n C_{\phi^n} \chi_{\{v\}}= \lambda^n \chi_{\phi^{-n}(\{v\})}$. Hence, if $|\lambda|<1$ we obtain that $(\lambda C_\phi)^n \chi_{\{v\}} \to 0$, as $n \to \infty$.  If $|\lambda|=1$, since $\{ n \in \N \, : \, \phi^{-n}(\{v\}) \neq \varnothing\}$ is finite we get that $(\lambda C_\phi)^n \chi_{\{v\}} = \lambda^n \chi_{\phi^{-n}(\{v\})} = 0$ for large enough $n$, so we also obtain that $(\lambda C_\phi)^n \chi_{\{v\}} \to 0$, as $n \to \infty$. Now, if $f \in X$, then $f$ is a finite linear combination of functions of the form $\chi_{\{v\}}$, and therefore it follows that $(\lambda C_\phi)^n f \to 0$, as $n \to \infty$, by the linearity of $(\lambda C_\phi)^n$.
\item Let $v \in T$ and and observe that $m(n,v) \geq 1$  since $\phi$ is injective. We define, for each $n \in \N$ and each function $\chi_{\{v\}}$, the function $B_n \chi_{\{v\}} \in \calF$ as
\[
  \hskip25pt (B_n \chi_{\{v\}})(u) =
  \begin{cases}
  \frac{m(n,v)-\dist(u,\phi^n(v))}{m(n,v)}, & \text{ if } 0 \leq \dist(u,\phi^n(v)) \leq m(n,v)-1,\\
  0,& \text{ otherwise,}
      \end{cases}
\]
and extend $B_n$ linearly to functions in $X$. Since $T$ is locally finite, the function $B_n \chi_{\{v\}}$ is of finite support and hence belongs to $X$; thus for each $f \in X$, we also have $B_n f \in X$.

  Now, choose $u \in T^*$ and let $j:=\dist(u,\phi^n(v))$. It then follows that $\dist(u^-,\phi^n(v))=j\pm 1$. If $\dist(u^-,\phi^n(v))=j+1$ and  $0 \leq j \leq m(n,v)-1$, this implies
\[
  \left|(B_n \chi_{\{ v\}})(u) -   (B_n \chi_{\{ v\}})(u^-)\right| = \frac{1}{m(n,v)}.
\]
If $\dist(u^-,\phi^n(v))=j-1$ and $0 < j \leq m(n,v)$, this implies
\[
  \left|(B_n \chi_{\{ v\}})(u) -   (B_n \chi_{\{ v\}})(u^-)\right| = \frac{1}{m(n,v)}.
\]
In any other case $|B_n \chi_{\{ v\}})(u) -   (B_n \chi_{\{ v\}})(u^-)|=0$. We obtain
\[
  \left|(B_n \chi_{\{ v\}})(u) -   (B_n \chi_{\{ v\}})(u^-)\right| \leq  \frac{1}{m(n,v)}, \quad \text{ for all } u \in T^*.
\]

Now, observe that, by hypothesis (2), we have $|\phi^n(v)| + c \geq m(n,v)$  for $n\geq N$ and hence $\dist(\root,\phi^n(v))=|\phi^n(v)| \geq m(n,v)-c$. It follows that
\[
  \frac{m(n,v)-\dist(\root,\phi^n(v))}{m(n,v)} \leq \frac{c}{m(n,v)},
  \]
  and hence by the definition of $B_n \chi_{\{ v \}}$ we have that
  \[
    0 \leq (B_n \chi_{\{ v \}})(\root)\leq \frac{c}{m(n,v)}, \quad \text{ for } n \geq N.
    \]
From this it follows that for $n \geq N$, we have
\begin{align*}
  \hskip35pt  \| B_n \chi_{\{v\}} \|_\Lip
  =& \max\Big\{|(B_n \chi_{\{v\}})(\root)|,  \sup_{u\in T^*} | B_n \chi_{\{v\}}(u) - B_n \chi_{\{v\}}(u^-)| \Big\} \\
  \leq&\frac{c}{m(n,v)},
\end{align*}
and hence
  \[
    \left\| \frac{1}{\lambda^n }B_n \chi_{\{v\}} \right\|_\Lip \leq \frac{c}{|\lambda|^n m(n,v)} \to 0 \quad\text{as } n\to \infty,
    \]
 by hypothesis (1). Hence, if $f \in X$, then
  \[
\frac{1}{\lambda^n }B_n f \to 0  \quad \text{ as } n \to \infty.
  \]
\item Now, for each $u \in T$, we have $C_{\phi^n} (B_n \chi_{\{v\}})(u)=(B_n \chi_{\{v\}})(\phi^n(u))$. By the definition of $B_n  \chi_{\{v\}}$ this is zero unless $0 \leq \dist(\phi^n(u),\phi^n(v))) \leq m(n,v)-1$. By definition of $m(n,v)$, this occurs only if $u=v$. But then $$(B_n \chi_{\{v\}})(\phi^n(u)) = (B_n \chi_{\{v\}})(\phi^n(v)) = 1$$ and hence
\[
    (\lambda C_\phi)^n \frac{1}{\lambda^n} B_n \chi_{\{v\}}
 = \chi_{\{v\}}.
\]

Hence $(\lambda C_\phi)^n \frac{1}{\lambda^n} B_n f =f$ for all $f \in X$.
\end{enumerate}
By the Hypercyclicity Criterion, it follows that $\lambda C_\phi$ is mixing.
\end{proof}

Observe that in the proof of this theorem, we did not use that $\phi$ has no periodic points. Nevertheless, this condition is not only necessary, but it follows from conditions (1) and (2) of Theorem~\ref{th:hyper_referee}.

Also, in Theorem~\ref{th:hyper_referee}, if conditions (1) and (2) hold only for the elements of a sequence $(n_k)$ of natural numbers, we obtain hypercyclicity instead of mixing.

From the above theorem, we obtain the following result.

\begin{example}\label{ex:example2}
 Consider the tree  $T=\N_0$, with $0$ being the root, and with two vertices $n$ and $m$ adjacent if and only if $|n-m|=1$. If $\phi:T \to T$ is defined by $\phi(m)=2m+1$, then $C_\phi$ is mixing.
\end{example}
\begin{proof}
  First of all, it is clear that for each $m \in \N_0$ the set $\{ n \in \N \, : \, \phi^{-n}(\{m\})\neq \varnothing
  \}$ is finite. As was shown in Example~\ref{ex:example1}, for any $n \in \N$ we have that $\phi^n(m)=2^n m + 2^n-1$. From this it follows that
  \[
    \min\{ \dist(\phi^n(k),\phi^n(m)) \, :\, k \in \N_0, k \neq m\}=2^n,
  \]
  and hence Condition (1) of Theorem~\ref{th:hyper_referee} holds. Also, $|\phi^n(m)|+1=2^n m + 2^n \geq 2^n$ for all $m \in \N_0$ and all $n \in \N$. Hence Condition (2) of Theorem~\ref{th:hyper_referee} also holds. It follows that $C_\phi$ is mixing.
\end{proof}

\subsection{Concluding remarks}

  Are there some conditions that are both sufficient and necessary for the hypercyclicity of $C_\phi$? We have not been able to solve this problem and we leave it open for future research. Another interesting question is whether the Lipschitz space and the little Lipschitz space can be defined if the tree is not locally finite, and what the consequences are for multiplication and composition operators on them. We plan to investigate this in the future.

\section{Acknowledgment}

We owe a debt of gratitude to the reviewers for their careful examination of the paper and their suggestions that greatly improved the paper. In particular, we wish to note the improvement of the proof of Theorem~\ref{th:Cphi_not} and the statement and proof of Theorems~\ref{th:CphiLip} and \ref{th:hyper_referee} suggested by anonymous reviewers.

\end{document}